\documentclass[letterpaper,11pt]{article}%
\usepackage{amsmath}
\usepackage{amsfonts}
\usepackage{graphicx}
\usepackage{amssymb}%
\newtheorem{theorem}{Theorem}[section]

\newtheorem{corollary}[theorem]{Corollary}

\newtheorem{definition}[theorem]{Definition}
\newtheorem{example}[theorem]{Example}

\newtheorem{lemma}[theorem]{Lemma}

\newtheorem{proposition}[theorem]{Proposition}
\newtheorem{remark}[theorem]{Remark}

\begin{document}

\begin{center}
{\Large A Note on Approximate Liftings}
\end{center}

\vspace{0.3cm}

\begin{center}
{Don Hadwin \hspace{3cm} Weihua Li}

\vspace{0.3cm}

Mathematics Department, University of New Hampshire, Durham, NH
03824

\vspace{0.3cm}

 {don@unh.edu \hspace{2cm} whli@unh.edu}
\end{center}

\noindent\textbf{Abstract:}\ In this paper, we prove approximate
lifting results in the C$^{\ast}$-algebra and von Neumann algebra
settings. In the C$^{\ast}$-algebra setting, we show that two
(weakly) semiprojective unital C*-algebras, each generated by $n$
projections, can be glued together with partial isometries to define
a larger (weakly) semiprojective algebra. In the von Neumann algebra
setting, we prove lifting theorems for trace-preserving
*-homomorphisms from abelian von Neumann algebras or hyperfinite von
Neumann algebras into ultraproducts. We also extend a classical
result of S. Sakai \cite{sakai} by showing that a tracial
ultraproduct of C*-algebras is a von Neumann algebra, which yields a
generalization of Lin's theorem \cite{Lin} on almost commuting
selfadjoint operators with respect to $\Vert\cdot\Vert_{p}$ on any
unital C*-algebra with trace.

\vspace{0.2cm}


\vspace{0.2cm}


\section{Introduction}

The key idea in this paper is the study of defining properties of
(or tuples of) operators such that an operator that
\textquotedblleft almost" satisfies this property is
\textquotedblleft close" to an operator that actually does satisfy
this property. In the setting of C*-algebras we would insist that
the \textquotedblleft closeness" be with respect to the norm, and in
the finite von Neumann algebra sense \textquotedblleft closeness"
would be with respect to the $2$-norm $\left\Vert \cdot\right\Vert
_{2}$ defined in terms of the
tracial state $\tau$ on the algebra, i.e., $\left\Vert x\right\Vert _{2}%
=\tau\left(  x^{\ast}x\right)  ^{1/2}$. A classic example of this
phenomenon is the fact that if $A$ is an operator such that
$\left\Vert A-A^{\ast }\right\Vert $ is small and $\left\Vert
A-A^{2}\right\Vert $ is small, then $A$ is very close to a
projection $P.$ In fact, $P$ can be chosen in the nonunital
C*-algebra generated by $A$. It is also true that if $\mathcal{M}$
is a finite von Neumann algebra with a faithful normal trace $\tau$
and if $A\in\mathcal{M}$ such that $\left\Vert A-A^{\ast}\right\Vert
_{2}$ is small and $\left\Vert A-A^{2}\right\Vert _{2}$ is small,
then there is a projection $P\in W^{\ast}\left(  A\right)  $ (the
von Neumann subalgebra generated by $A$ in $\cal M$) such that
$\left\Vert A-P\right\Vert _{2}$ is small.

In the C*-algebra setting these ideas are essentially the notion of weak
semiprojectivity introduced by S. Eilers and T. Loring \cite{Ei} in 1999, and
semiprojectivity introduced by B. Blackadar \cite{B} in 1985. These notions
were studied by T. Loring \cite{Lo1} in terms of stable relations and D.
Hadwin, L. Kaonga and B. Mathes \cite{Don} in terms of their noncommutative
continuous functions.

The von Neumann algebra results appear in an ad hoc manner in various papers
in the literature.

Semiprojectivity and weak semiprojectivity can also be expressed in
terms of liftings of representations into algebras of the form
$\prod_{1}^{\infty}{\mathcal{B}}_{n}/\oplus
_{1}^{\infty}{\mathcal{B}}_{n}$ or in terms of ultraproducts $\prod
_{1}^{\infty}{\mathcal{B}}_{n}/\mathcal{J}$. It is in the theory of
(tracial) ultraproducts of finite von Neumann algebras where many of
these ``approximate" results appear in the von Neumann algebra
setting.

After the preliminary definitions and results (Section 2), we begin Section 3
with our results in the C*-algebra setting. Our main result is that two
(weakly) semiprojective unital C*-algebras, each generated by $n$ projections,
can be glued together with partial isometries to define a larger (weakly)
semiprojective algebra (Theorem \ref{theorem,glue together, projections}).

In the von Neumann algebra setting (Section 4) we prove lifting
theorems for trace-preserving *-homomorphisms from abelian von
Neumann algebras (Corollary \ref{corollary,tracepreserving
homomorphism}) or hyperfinite von Neumann algebras (Theorem
\ref{theorem,hyperfinite}) into ultraproducts. We also extend a
classical result of S. Sakai \cite{sakai} by showing (Theorem
\ref{theorem, saikai}) that a tracial ultraproduct of C*-algebras is
a von Neumann algebra. This result allows us to prove a hybrid
result (Corollary \ref{corollary,shanghai}), namely, an approximate
result with respect to $\left\Vert \cdot\right\Vert _{2}$ on
C*-algebras. For example, if $\varepsilon>0,$ then there is a
$\delta>0$ such that for any unital C*-algebra $\mathcal{A}$ with
trace $\tau,$ we have that if $u,v$ are unitaries in $\mathcal{A}$
with $\left\Vert uv-vu\right\Vert _{2}<\delta$, then there are
commuting unitaries $u^{\prime},v^{\prime}$ in C$^{\ast}\left(
u,v\right)  $ (the unital C$^*$-subalgebra generated by $u, v$ in
$\cal A$) such that $\left\Vert u-u^{\prime}\right\Vert
_{2}+\left\Vert v-v^{\prime }\right\Vert _{2}<\varepsilon$. With
respect to the operator norm this fails even in the class of
finite-dimensional C*-algebras \cite{Vo}.

\section{Preliminaries}

A C*-algebra $\mathcal{A}$ is \emph{projective} if, for any
*-homomorphism $\varphi:\mathcal{A}\rightarrow\mathcal{C}$, where
$\mathcal{C}$ is a C$^{*}$-algebra, and every surjective
*-homomorphism $\rho:\mathcal{B}\rightarrow\mathcal{C}$, where $\cal C$ is a C$^*$-algebra, there is a
*-homomorphism $\overline{\varphi}:\mathcal{A}\rightarrow\mathcal{B}$
such that $\rho \circ\overline{\varphi}=\varphi$. A
C$^{\ast}$-algebra $\mathcal{A}$ is \emph{semiprojective} \cite{B}
if, for every *-homomorphism $\pi:\
{\mathcal{A}}\rightarrow{\mathcal{B}}/\overline{\cup_{1}^{\infty
}{\mathcal{J}}_{n}}$, where ${\mathcal{J}}_{n}$ are increasing
ideals of a
C$^{\ast}$-algebra $\mathcal{B}$, and with $\varphi_{N}:\ {\mathcal{B}%
}/{\mathcal{J}}_{N}\rightarrow{\mathcal{B}}/\overline{\cup_{1}^{\infty
}{\mathcal{J}}_{n}} $ the natural quotient map, there exists a $\ast
$-homomorphism $\pi_{N}:\ {\mathcal{A}}\rightarrow{\mathcal{B}}/{\mathcal{J}%
}_{N}$ such that $\pi=\varphi_{N}\circ\pi_{N}$. A C$^{\ast}$-algebra
$\mathcal{A}$ is \emph{weakly semiprojective} \cite{Lo1} if, for any given
sequence $\left\{  \mathcal{B}_{n}\right\}  _{n\in\mathbb{N}}$ of C$^{\ast}%
$-algebras and a *-homomorphism $\pi:\ {\mathcal{A}}\rightarrow\prod
_{1}^{\infty}{\mathcal{B}}_{n}/\oplus_{1}^{\infty}{\mathcal{B}}_{n}$,
there exist functions
$\pi_{n}:\mathcal{A}\rightarrow\mathcal{B}_{n}$ for all $n\geq1$ and
a positive integer $N$ such that

(1) \ $\pi_{n}$ is a *-homomorphism for all $n\geq N,$ and

(2)\ $\pi\left(  a\right)  =[\left\{  \pi_{n}\left(  a\right)
\right\}]$ for every $a\in\mathcal{A}$.
\newline Equivalently, since
$\prod_{1}^{\infty}{\mathcal{B}}_{n}/\oplus_{1}^{\infty
}{\mathcal{B}}_{n}$ is isomorphic to $\prod_{N}^{\infty}{\mathcal{B}}%
_{n}/\oplus_{N}^{\infty}{\mathcal{B}}_{n},$ the conditions above say
that there is a *-homomorphism
$\rho:\mathcal{A}\rightarrow\prod_{N}^{\infty }{\mathcal{B}}_{n}$
such that $\pi\left(  a\right)  =\rho\left(  a\right)
+\oplus_{N}^{\infty}{\mathcal{B}}_{n}$ for every $a\in\mathcal{A}$.

These notions of projectivity makes sense in two categories:

(1)\ the \emph{nonunital category}, i.e., the category of C*-algebras with
*-homomorphisms as morphisms, and

(2)\ the \emph{unital category}, i.e., the category of unital
C*-algebras with unital *-homomorphisms as morphisms.

These notions are drastically different in the different categories. For
example, the $1$-dimensional C*-algebra $\mathbb{C}$ is projective in the
unital category, but not in the nonunital category, e.g., in the definition of
projective C$^{*}$-algebra, let $\mathcal{B}=C_{0}\left(  (0,1]\right)  ,$
$\mathcal{C=}\mathbb{C}$ and $\rho\left(  f\right)  =f\left(  1\right)  $.
However, if $\mathcal{A}$ is not unital and projective (semiprojective, weakly
semiprojective) in the nonunital category, and if $\mathcal{A}^{+}$ is the
algebra obtained by adding a unit to $\mathcal{A},$ then $\mathcal{A}^{+}$ is
projective (semiprojective, weakly semiprojective) in the unital category. In
Loring's book \cite{Lo1} he only considers the nonunital category. In this
paper we restrict ourselves to the unital category.

Suppose $\mathcal{S}$ is a subset of a unital C*-algebra $\cal A$.
Let C$^{\ast}\left( \mathcal{S}\right)  $ denote the unital
C*-subalgebra generated by $\mathcal{S}$ in $\cal A$.

In \cite{Don} the notions of semiprojectivity and weak
semiprojectivity for finitely generated algebras were cast in terms
of \emph{noncommutative continuous functions}. The *-algebra of
noncommutative continuous functions is basically the metric
completion of the algebra of $\ast $-polynomials with respect to a
family of seminorms. There is a functional calculus for these
functions on any $n$-tuple of operators on any Hilbert space. Here
is a list of a few of the basic properties of noncommutative
continuous functions \cite{Don}:

(1)\ For each noncommutative continuous function $\varphi$ there is
a sequence $\left\{  p_{k}\right\}  $ of noncommutative
*-polynomials such that for every tuple $\left(
T_{1},\ldots,T_{n}\right)  $ we have
\[
\left\Vert p_{k}\left(  T_{1},\ldots,T_{n}\right)  -\varphi\left(
T_{1},\ldots,T_{n}\right)  \right\Vert \rightarrow0,
\]
and the convergence is uniform on bounded $n$-tuples of operators.

(2)\ For any tuple $\left(  T_{1},\ldots,T_{n}\right)  ,$ C*$\left(
T_{1},\ldots,T_{n}\right)  $ is the set of all $\varphi\left(  T_{1}%
,\ldots,T_{n}\right)  $ with $\varphi$ a noncommutative continuous function.

(3)\ For any $n$-tuple $\left(  A_{1},\ldots,A_{n}\right)  $ and any $S\in
C^{\ast}\left(  A_{1},\ldots,A_{n}\right)  ,$ there is a noncommutative
continuous function $\varphi$ such that $S=\varphi\left(  A_{1},\ldots
,A_{n}\right)  $ and $\left\Vert \varphi\left(  T_{1},\ldots,T_{n}\right)
\right\Vert \leq\left\Vert S\right\Vert $ for all $n$-tuples $\left(
T_{1},\ldots,T_{n}\right)  .$

(4)\ If $T_{1},\ldots,T_{n}$ are elements of a unital C*-algebra
$\mathcal{A} $ and $\pi:\mathcal{A}\rightarrow\mathcal{B}$ is a
unital *-homomorphism, then
\[
\pi\left(  \varphi\left(  T_{1},\ldots,T_{n}\right)  \right)  =\varphi\left(
\pi\left(  T_{1}\right)  ,\ldots,\pi\left(  T_{n}\right)  \right)
\]
for every noncommutative continuous function.

In \cite{Don} it was shown that the natural notion of relations used to define
a C*-algebra generated by $a_{1},\ldots,a_{n}$ are all of the form%
\[
\varphi\left(  a_{1},\ldots,a_{n}\right)  =0
\]
for a noncommutative continuous function $\varphi.$ In fact, it was also shown
in \cite{Don} that given a unital C*-algebra $\mathcal{A}$ generated by
$a_{1},\ldots,a_{n},$ there is a single noncommutative continuous function
$\varphi$ such that $\mathcal{A}$ is isomorphic to the universal C*-algebra
C*$\left(  x_{1},\ldots,x_{n}|\varphi\right)  $ with generators $x_{1}%
,\ldots,x_{n}$ and with the single relation $\varphi\left(  x_{1},\ldots
,x_{n}\right)  =0,$ where the map $x_{j}\mapsto a_{j}$ extends to a $\ast
$-isomorphism. Such a noncommutative continuous function $\varphi$ must be
\emph{null-bounded}, i.e., there is a number $r>0$ such that $\left\Vert
A_{j}\right\Vert \leq r$ for $1\leq j\leq n$ whenever $\varphi\left(
A_{1},\ldots,A_{n}\right)  =0$. In this sense, every finitely generated
C*-algebra is finitely presented. In particular, Theorem 14.1.4 in T. Loring's
book \cite{Lo1} is true for all finitely generated C*-algebras.

For a finitely generated nonunital C*-algebra $\mathcal{A}$ there is a
null-bounded noncommutative continuous function $\varphi$ such that
$\mathcal{A}$ is isomorphic to the universal nonunital C*-algebra C$_{0}%
^{\ast}\left(  x_{1},\ldots,x_{n}|\varphi\right)$ with generators
$x_{1},\ldots,x_{n}$ and with the single relation $\varphi\left(  x_{1}%
,\ldots,x_{n}\right)  =0$.

Here is a reformulation of the notions of semiprojectivity and weak
semiprojectivity for finitely generated C*-algebras in terms of noncommutative
continuous functions. We only state the result in the unital category.

\begin{proposition}
\cite{Don} Suppose $\varphi$ is a null-bounded noncommutative continuous
function. Then

(1)\ C*$\left(  x_{1},\ldots,x_{n}|\varphi\right)  $ is weakly semiprojective
if and only if there exist noncommutative continuous functions $\varphi
_{1},\ldots,\varphi_{n}$ such that for any $\varepsilon>0$, there exists
$\delta>0$, such that for any operators $T_{1},\cdots,T_{n}$ with
$\Vert\varphi(T_{1},\ldots,T_{n})\Vert<\delta$, we have

\hspace{2em}(a) $\varphi(\varphi_{1}(T_{1},\ldots,T_{n}),\ldots,\varphi
_{n}(T_{1},\ldots,T_{n}))=0$, and

\hspace{2em}(b) $\Vert T_{j}-\varphi_{j}(T_{1},\ldots,T_{n})\Vert<\varepsilon$.

(2)\ C*$\left(  x_{1},\ldots,x_{n}|\varphi\right)  $ is semiprojective if, in
addition, we can choose $\varphi_{1},\ldots,\varphi_{n}$ as in part 1 so that
$\varphi_{j}\left(  A_{1},\ldots,A_{n}\right)  =A_{j}$ for $1\leq j\leq n$,
whenever $\varphi\left(  A_{1},\ldots,A_{n}\right)  =0$.
\end{proposition}

We call the functions $\varphi_{1},\ldots,\varphi_{n}$ the \emph{(weakly)
semiprojective approximating functions} for $\varphi$.

For example, it is a classical result that has often been rediscovered that a
selfadjoint operator $A$ with $\left\Vert A-A^{2}\right\Vert $ sufficiently
small is very close to a projection. More precisely, if $\left\Vert
A-A^{2}\right\Vert <\varepsilon^{2}<1/9,$ then $\sigma\left(  A\right)
\subset\left(  -\varepsilon,\varepsilon\right)  \cup\left(  1-\varepsilon
,1+\varepsilon\right)  $. So if $h:\mathbb{R}\rightarrow\mathbb{R}$ is defined
by
\[
h\left(  t\right)  =\left\{
\begin{array}
[c]{ll}%
0, & \mbox{if }t\leq\frac{1}{3}\\
3t-1, & \mbox{if }\frac{1}{3}<t<\frac{2}{3}\\
1, & \mbox{if }\frac{2}{3}\leq t
\end{array}
\right.  ,
\]
then $h\left(  A\right)  $ is a projection and $\left\Vert A-h\left(
A\right)  \right\Vert =\sup\limits_{t\in\sigma\left(  A\right)
}\left\vert t-h\left(  t\right)  \right\vert <\varepsilon$. If $A$
already is a projection, then $h\left(  A\right)  =A.$ Thus the
universal C*-algebra generated by a single projection is C*$\left(
x|\varphi\right)  $ where $\varphi\left(  x\right)  =\left(
x-x^{\ast}\right)  ^{2}+\left( x-x^{2}\right)  ^{\ast}\left(
x-x^{2}\right)  ,$ and defining $\varphi _{1}\left(  x\right)
=h\left( \frac{x+x^*}{2}  \right)  $ shows that C*$\left(
x|\varphi\right) $ is semiprojective.

Throughout this paper, all the C$^*$-algebras considered are unital
and all *-homomorphism are unital.

\section{C*-algebra Results}

For simplicity we only consider finitely generated C$^{*}$-algebras throughout
this section.

The main results in this section concern the (weak) semiprojectivity of
C*-algebras defined in terms of partial isometries. We begin with some results
that are elementary in the unital category.

\begin{lemma}
\label{lemma, basic properties} Suppose $\mathcal{A}$ and $\mathcal{B}$ are
separable unital C*-algebras. The following are true:

(1)\ if $\mathcal{A}$ is (weakly) semiprojective, then $\mathcal{A}%
\otimes\mathcal{M}_{n}\left(  \mathbb{C}\right)  $ is (weakly) semiprojective;

(2)\ if $\mathcal{A}\otimes\mathcal{M}_{n}\left(  \mathbb{C}\right)  $ is
weakly semiprojective, then $\mathcal{A}$ is weakly semiprojective;

(3)\ if $\mathcal{A}$ and $\mathcal{B}$ are projective (semiprojective, weakly
semiprojective), then so is $\mathcal{A}\ast\mathcal{B}$;

(4)\ if $\mathcal{A}\ast\mathcal{B}$ is (weakly semiprojective,
semiprojective) projective, and there is a linear multiplicative functional
$\alpha$ on $\mathcal{B}$, then $\mathcal{A}$ is (weakly semiprojective,
semiprojective) projective;

(5)\ $\mathcal{A}\oplus\mathcal{B}$ is (weakly) semiprojective if an only if
both $\mathcal{A}$ and $\mathcal{B}$ are (weakly) semiprojective.
\end{lemma}

Proof. \ (1)\ Suppose $\mathcal{A}$ is weakly semiprojective, and
$\mathcal{A}=C^{\ast}(x_{1},\ldots,x_{m}|\varphi)$ with weakly semiprojective
approximating functions $\varphi_{1},\ldots,\varphi_{m}$. Since $\mathcal{M}%
_{n}(\mathbb{C})$ is semiprojective, we can assume that ${\cal M}_{n}%
(\mathbb{C})=C^{\ast}(y|\rho)$ with a semiprojective approximating function
$\rho_{1}$. Hence
\[
\mathcal{A}\otimes\mathcal{M}_{n}\left(  \mathbb{C}\right)  =C^{\ast}%
(x_{1},\ldots,x_{m},y|\Phi)
\]
where
\begin{align*}
\Phi(x_{1},\ldots,x_{m},y) & =\sum_{i=1}^{m}(x_{i}y-yx_{i})^{\ast}%
(x_{i}y-yx_{i})+\sum_{i=1}^{m}(x_{i}y^{\ast}-y^{\ast}x_{i})^{\ast}%
(x_{i}y^{\ast}-y^{\ast}x_{i})\\
& +\varphi(x_{1},\ldots,x_{m})^{\ast}\varphi(x_{1},\ldots,x_{m})+\rho
(y)^{\ast}\rho(y).
\end{align*}

Since matrix units $E_{i,j}(1\leq i,j\leq n)$ are in ${\mathcal{M}}%
_{n}(\mathbb{C})$, there exists a family of noncommutative
continuous functions $\{\rho_{i,j}: 1\leq i,j\leq n\}$ such that
$E_{i,j}=\rho_{i,j}(y)$.

For any operators $T_{1}, \ldots, T_{m}, S$, and for $1\leq j\leq m$, let
$\widehat{T_{k}}=\sum_{j=1}^{n}\rho_{j,1}(\rho_{1}(S))\cdot T_{k}\cdot
\rho_{1,j}(\rho_{1}(S))$. Define functions $\{\Phi_{k}: 1\leq k\leq m+1\}$ by
\[
\Phi_{k}\left( T_{1}, \ldots, T_{m}, S\right)  =\left\{
\begin{array}
[c]{ll}%
\varphi_{k}\left( \widehat{T_{1}}, \ldots, \widehat{T_{m}} \right)  , & 1\leq
k\leq m\\
\rho_{1}(S), & k=m+1
\end{array}
\right.
\]

Given any $\varepsilon>0$, we will find some $\delta>0$ in the definition of
weak semiprojectivity.

Note that $\mathcal{M}_{n}(\mathbb{C)}$ is semiprojective, there exists
$\delta_{1}>0$, such that if $\|\rho(S)\|<\delta_{1}$, then

(1)\ $\rho(\Phi_{m+1}(T_{1}, \ldots, T_{m}, S))=\rho(\rho_{1}(S))=0$,

(2)\ $\|\rho_{1}(S)-S\|<\varepsilon$,

(3)\ $\rho_{1}(S)$ and $\rho_{1}(S)^{*}$ commute with all $\widehat{T_{k}}$.

Since $\mathcal{A}$ is weakly semiprojective, there exists $\delta_{2}>0 $,
such that if $\|\varphi(\widehat{T_{1}},\cdots, \widehat{T_{m}})\|<\delta_{2}
$, then
\[
\varphi(\varphi_{1}(\widehat{T_{1}},\cdots, \widehat{T_{m}}), \cdots,
\varphi_{m}(\widehat{T_{1}},\cdots, \widehat{T_{m}}))=0
\]
and
\[
\|\widehat{T_{k}}-\varphi_{k}(\widehat{T_{1}},\cdots, \widehat{T_{m}%
})\|<\varepsilon.
\]
Note that
\[
\|\varphi(\widehat{T_{1}},\cdots, \widehat{T_{m}})\|\leq\|\varphi
(\widehat{T_{1}},\cdots, \widehat{T_{m}})-\varphi(T_{1}, \ldots,
T_{m})\|+\|\varphi(T_{1}, \ldots, T_{m})\|,
\]
there exists $\delta_{3}>0$, such that if $\|T_{k}-\widehat{T_{k}}%
\|<\delta_{3}$, then
\[
\|\varphi(\widehat{T_{1}},\cdots, \widehat{T_{m}})\|\leq\|\varphi
(\widehat{T_{1}},\cdots, \widehat{T_{m}})-\varphi(T_{1}, \ldots,
T_{m})\|<\frac{\delta_{2}}{2}.
\]
In addition, if $\|\varphi(T_{1}, \ldots, T_{m})\|<\frac{\delta_{2}}{2} $,
then $\|\varphi(\widehat{T_{1}},\cdots, \widehat{T_{m}})\|<\delta_{2}$.
Furthermore,
\begin{align*}
\|T_{k}-\widehat{T_{k}}\|  &  =\|T_{k}-\sum_{j=1}^{n}\rho_{j,1}(\rho
_{1}(S))\cdot T_{k}\cdot\rho_{1,j}(\rho_{1}(S))\|\\
&  =\|\sum_{j=1}^{n}\rho_{j,1}(\rho_{1}(S))\cdot\left(  \rho_{1,j}(\rho
_{1}(S))\cdot T_{k}-T_{k}\cdot\rho_{1,j}(\rho_{1}(S))\right)  \|\\
&  \leq\sum_{j=1}^{n}\|\rho_{1,j}(\rho_{1}(S))\cdot T_{k}-T_{k}\cdot\rho
_{1,j}(\rho_{1}(S))\|.
\end{align*}
Therefore there exists $\delta_{4}>0$, such that if $\|\rho_{1,j}(\rho
_{1}(S))\cdot T_{k}-T_{k}\cdot\rho_{1,j}(\rho_{1}(S))\|<\delta_{4}$, then
$\|T_{k} -\widehat{T_{k}}\|<\delta_{3}$.

Note that $\rho_{1}(S)=\sum_{i,j=1}^{n}c_{ij}\cdot\rho_{ij}(\rho_{1}(S))$,
where $c_{i,j}$'s are complex numbers. We have
\begin{align*}
\|T_{k}\rho_{1}(S)-\rho_{1}(S)T_{k}\|  &  =\|T_{k}\sum_{i,j=1}^{n}
c_{i,j}\cdot\rho_{ij}(\rho_{1}(S))-\sum_{i,j=1}^{n}c_{i,j}\cdot\rho_{ij}%
(\rho_{1}(S))T_{k}\|\\
&  \leq\sum_{i,j=1}^{n}|c_{ij}|\cdot\|T_{k} \rho_{ij}(\rho_{1}(S))-\rho
_{ij}(\rho_{1}(S))T_{k}\|\\
&  =\sum_{i,j=1}^{n}|c_{ij}|\cdot\delta_{4}.
\end{align*}

Let $\delta_{5}=\sum_{i,j=1}^{n}|c_{ij}|\cdot\delta_{4}$. Since
\begin{align*}
\|T_{k}\rho_{1}(S)-\rho_{1}(S)T_{k}\|  &  =\|T_{k}\left(  \rho_{1}%
(S)-S+S\right)  -\left(  \rho_{1}(S)-S+S\right)  T_{k}\|\\
&  \leq2\|T_{k}\|\cdot\|\rho_{1}(S)-S\|+\|T_{k}S-ST_{k}\|,
\end{align*}
there exists $\delta_{6}>0$, such that if $\|\rho_{1}(S)-S\|<\delta_{6},$
$\|T_{k}S-ST_{k}\|<\delta_{6}$, then $\|T_{k}\rho_{1}(S)-\rho_{1}%
(S)T_{k}\|<\delta_{5}$.

By the fact that ${\cal M}_n(\Bbb C)=C^*(y|\rho)$ is semiprojective,
there exists $\delta_{7}>0$ such that if
$\|\varphi(S)\|<\delta_{7}$, then $\|\rho_{1}(S)-S\|<\delta_{6}$.

Note that $\Vert\varphi\left(  T_{1},\ldots,T_{m}\right)  \Vert,\Vert
\rho(S)\Vert,\Vert T_{i}S-ST_{i}\Vert,\Vert T_{i}S^{\ast}-S^{\ast}T_{i}\Vert$
are all less than or equal to $\sqrt{\Vert\Phi\left(  T_{1},\ldots
,T_{m},S\right)  \Vert}$. Put $\delta=\mbox{min}\{\delta_{1}^{2},(\delta
_{2}/2)^{2},\delta_{6}^{2},\delta_{7}^{2}\}$, then $\Phi_{1},\ldots,\Phi
_{m+1}$ are weakly semiprojective approximating functions for $\Phi$.

If $\mathcal{A}$ is semiporjective and $\varphi_{1},\ldots,\varphi_{m}$ are
semiprojective approximating functions for $\varphi$, it is clear that
$\Phi_{1},\ldots,\Phi_{m+1}$ are semiprojective approximating functions for
$\Phi$.

(2)\ Suppose ${\mathcal{A}}\otimes{\mathcal{M}}_{n}(\mathbb{C)}$ is
weakly
semiprojective. Let $\pi:\mathcal{A}\rightarrow\prod_{1}^{\infty}{\mathcal{B}%
}_{k}/\oplus_{1}^{\infty}{\mathcal{B}}_{k}$ is a unital
*-homomorphism. Then $\rho=\pi\otimes id$ is a unital *-homomorphism from $\mathcal{A}\otimes{\mathcal{M}}_{n}(\mathbb{C)}$
 to $\left(  \prod_{1}^{\infty}{\mathcal{B}}_{k}/\oplus_{1}^{\infty
}{\mathcal{B}}_{k}\right) \otimes{\mathcal{M}}_{n}(\mathbb{C)}\
=\prod _{1}^{\infty}\left(
{\mathcal{B}}_{k}\otimes{\mathcal{M}}_{n}(\mathbb{C)} \right)
/\oplus_{1}^{\infty}\left(  {\mathcal{B}}_{k}\otimes{\mathcal{M}}
_{n}(\mathbb{C)}\right)$. Since $\mathcal{A}\otimes{\mathcal{M}}
_{n}(\mathbb{C)}$ is weakly semiprojective, there is a positive
integer $N$ and maps
$\rho_{k}:\mathcal{A}\otimes{\mathcal{M}}_{n}(\mathbb{C)}\rightarrow
{\mathcal{B}}_{k} $ such that, for $k\geq N,$ $\rho_{k}$ is a unital
*-homomorphism and, for every $x\in\mathcal{A}\otimes{\mathcal{M}%
}_{n}(\mathbb{C)},$
\[
\rho\left(  x\right)  =\left[  \left\{  \rho_{k}\left(  x\right)  \right\}
\right]  .
\]
It follows that there is a sequence $\left\{  U_{k}\right\}  $ of
unitary elements ($U_{k}\in\mathcal{B}_{k}$) such that $\left\Vert
U_{k}-1\right\Vert \rightarrow0$ and $\left\Vert 1\otimes
T-U_{k}^{\ast}\rho_{k}\left(  1\otimes T\right)  U_{k}\right\Vert
\rightarrow0$ for every $T\in\mathcal{M}_{n}\left( \mathbb{C}\right)
$. Therefore, for $k\geq N$ and $A\in\mathcal{A}$,
$U_{k}^{\ast}\rho_{k}\left(  A\otimes1\right)  U_{k}$ is in the
commutant
of $1\otimes{\mathcal{M}}_{n}(\mathbb{C)},$ which is $\mathcal{B}_{k}%
\otimes1.$ Hence there are representations $\pi_{k}:\mathcal{A}\rightarrow
\mathcal{B}_{k}$ such that $\pi_{k}\left(  A\right)  \otimes1=U_{k}^{\ast}%
\rho_{k}\left(  A\otimes1\right)  U_{k}$ for every $A\in\mathcal{A}$. Clearly,
$\pi\left(  A\right)  =\left[  \left\{  \pi_{k}\left(  A\right)  \right\}
\right]  $ for every $A\in\mathcal{A}$.

(3)\ This is obvious from the defining properties of the free
product in the unital category.

(4) We give a proof for the projective case; the other cases are handled
similarly. Suppose $\mathcal{C}$ is a unital C*-algebra with an ideal
$\mathcal{J}$ and $\pi:\mathcal{A}\rightarrow\mathcal{C}/\mathcal{J}$ is a
*-homomorphism. Define a unital *-homomorphism $\sigma:\mathcal{B}%
\rightarrow\mathcal{C}/\mathcal{J}$ by $\sigma\left(  x\right)  =\alpha\left(
x\right)\cdot 1 $. Thus there is a unital *-homomorphism $\rho:\mathcal{A}%
\ast\mathcal{B}\rightarrow\mathcal{C}/\mathcal{J}$ such that $\rho
|\mathcal{A}=\pi$ and $\rho|\mathcal{B}=\sigma$. Since $\mathcal{A}%
\ast\mathcal{B}$ is projective, $\rho$ lifts to a *-homomorphism
$\tau:\mathcal{A}\ast\mathcal{B}\rightarrow\mathcal{C}$. Thus $\tau
|\mathcal{A}$ is the required lifting of $\pi$.

(5) Suppose $\mathcal{A}=C^{*}(x_{1}, \ldots, x_{m}|\varphi)$ is
weakly semiprojective with weakly semiprojective approximating
functions $\varphi_{1}, \ldots, \varphi_{m}$, and
$\mathcal{B}=C^{*}(y_{1}, \ldots, y_{n}|\psi)$ is weakly
semiprojective with weakly semiprojective approximating functions
$\psi_{1}, \ldots, \psi_{n}$. Then
\[
{\mathcal{A}}\oplus{\ \mathcal{B}}=C^{\ast}(x_{1},\ldots,x_{m},y_{1}%
,\ldots,y_{n},p|\Phi),
\]
where
\begin{align*}
& \Phi(x_{1},\ldots,x_{m},y_{1},\ldots,y_{n},p)\\
 =&\varphi(x_{1},\ldots,x_{m})^{\ast}\varphi(x_{1},\ldots,x_{m})+\psi
(y_{1},\ldots,y_{n})^{\ast}\psi(y_{1},\ldots,y_{n})+\\
& +(p-p^{\ast})^{\ast}(p-p^{\ast})+(p-p^{2})^{\ast}(p-p^{2})+\sum_{j=1}%
^{m}(px_{j}-x_{j}p)^{\ast}(px_{j}-x_{j}p)+\\
& +\sum_{j=1}^{m}(px_{j}-x_{j})^{\ast}(px_{j}-x_{j})+\sum_{j=1}^{n}%
(py_{j}-y_{j}p)^{\ast}(py_{j}-y_{j}p).
\end{align*}

Let $f$ be a continuous function on $\mathbb{R}$ such that $f(t)=0$ when
$t\leq\frac{1}{4}$ and $f(t)=1$ when $t\geq\frac{3}{4}$. Define $\Phi_{1},
\ldots, \Phi_{m+n+1}$ on ${\mathcal{A}}\oplus{\mathcal{B}}$ by
\[
\Phi_{i}(x_{1}, \ldots, x_{m}, y_{1}, \ldots, y_{n},p)=\left\{
\begin{array}
[c]{ll}%
f(p)\varphi_{i}(x_{1}, \ldots, x_{m})f(p), & 1\leq i\leq m\\
(1-f(p))\psi_{i-m}(y_{1}, \ldots, y_{n})(1-f(p)), & m+1\leq i\leq m+n\\
f(p), & i=m+n+1
\end{array}
\right.
\]

For any $\varepsilon>0$, and any operators $S_{1}, \ldots,
S_{m},T_{1}, \ldots, T_{n}, Q $, there exists $\delta_{1}>0$, such
that if $\|P-P^{*}\|<\delta_1$ and $\|Q-Q^{2}\|<\delta_{1}$, then
$f(Q) $ is a projection and $\|Q-f(Q)\|<\varepsilon$. In addition,
there exists $\delta_{2}>0$, such that if $\|\varphi(S_{1}, \ldots,
S_{m})\|<\delta_{2}$, then $$\|S_{j}-\varphi _{j}(S_{1}, \ldots,
S_{m} )\|<\varepsilon\ \ \ \mbox{ for all}\ \ \  1\leq j\leq m$$ and
$$\varphi\left( \varphi_{1}(S_{1}, \ldots, S_{m}), \ldots,
\varphi_{m}(S_{1}, \ldots, S_{m})\right)  =0;$$ there exists
$\delta_{3}>0$, such that if $\|\psi(T_{1}, \ldots,
T_{n})\|<\delta_{3}$, then $$\|T_{k}-\psi_{k}(T_{1}, \ldots,
T_{n})\|<\varepsilon\ \ \ \mbox{ for all} \ \ \ 1\leq k\leq n$$ and
$$\psi\left( \psi_{1}(T_{1}, \ldots, T_{n}), \ldots, \psi_{n}(T_{1},
\ldots, T_{n}) \right)  =0.$$ Let
$$\widehat{S_{j}}=f(Q)\varphi_{j}(S_{1}, \ldots, S_{m})f(Q)\ \ \mbox{ for
all}\ \ 1\leq j\leq m$$ and
$$\widehat{T_{k}}=(1-f(Q))\psi_{k}(T_{1}, \ldots,
T_{n})(1-f(Q))\ \ \mbox{ for all} \ \ 1\leq k\leq n,$$ then $\widehat{S_{j}}%
f(Q)=f(Q)\widehat{S_{j}}=\widehat{S_{j}}$ and $\widehat{T_{k}}%
f(Q)=f(Q)\widehat{T_{k}}=0$.

Choose
$\delta=\mbox{min}\{\delta_{1}^{2},\delta_{2}^{2},\delta_{3}^{2}\}$,
then $\Phi_{1}, \ldots, \Phi_{m+n+1}$ are weakly semiprojective
approximating functions for $\Phi$.

Conversely, suppose $\mathcal{A}\oplus\mathcal{B}=C^{\ast}(x_{1}\oplus
y_{1},\ldots,x_{n}\oplus y_{n}|\varphi)$ is weakly semiprojective with weakly
semiprojective approximating functions $\varphi_{1}, \ldots, \varphi_{n}$.
Since $\varphi(A\oplus B)=\varphi(A)\oplus\varphi(B)$, it is clear that
${\mathcal{A}}=C^{\ast}(x_{1}, \ldots, x_{n}|\varphi)$ and ${\mathcal{B}%
}=C^{\ast}(y_{1}, \ldots, y_{n}|\varphi)$, and $\mathcal{A} $ and
$\mathcal{B}$ are all weakly semiprojective with weakly semiprojective
approximating functions $\varphi_{1}, \ldots, \varphi_{n}$.

It is not hard to prove the semiprojective case using the similar
idea.\hfill$\Box$

\begin{remark}
(1)\ Statement $\left(  2\right)  $ and statement $\left(  5\right)
$ of Lemma \ref{lemma, basic properties} are not true in the
nonunital case, and statement $\left(  1\right)  $ is not true for
projectivity in the unital case, e.g., in the Calkin algebra the
C$^{*}$-algebra generated by $\left(
\begin{array}
[c]{ll}%
0 & S\\
0 & 0
\end{array}
\right) $, where $S$ is the unilateral shift operator, is isomorphic to
$\mathcal{M}_{2}(\mathbb{C})$, but it cannot be lifted to a representation of
$\mathcal{M}_{2}(\mathbb{C})$ in $\mathcal{B(H)}$.

(2)\ In general the different types of projectivity are not
preserved under tensor products even when the algebras are very
nice. For example, if $X$ is the unit circle, then $C\left(
X\right)$, which is isomorphic to the universal C$^*$-algebra
generated by a unitary operator, is projective, but $C\left(
X\right) \otimes C\left( X\right) =C^{\ast }\left( x,y|x,y\text{
unitary, }xy-yx=0\right)  $ is not weakly semiprojective \cite{Vo}.

(3)\ In the nonunital category, statement $\left(  4\right)  $ of Lemma
\ref{lemma, basic properties} always holds, because the $0$ functional is allowed.
\end{remark}

The following lemma is a key ingredient to our main results in this section.

\begin{lemma}
\label{lemma,for theorem in section3} There exists a noncommutative continuous
function $\psi(x,y,z)$ such that, for any C*-algebra $\mathcal{A} $ and any
$\varepsilon>0$, there exists $\delta>0$, such that, whenever $P,Q,A\in
\mathcal{A}$ with $P$ and $Q$ projections and $\Vert A^{\ast}A-P\Vert<\delta$,
$\Vert AA^{\ast}-Q\Vert<\delta$, we have

(1)\ $\Vert\psi(P,Q,A)-A\Vert<\varepsilon$,

(2)\ $\psi(P,Q,A)^{\ast}\psi(P,Q,A)=P$ and $\psi(P,Q,A)\psi(P,Q,A)^{\ast}=Q$,

(3)\ $\psi(P,Q,A)=A$ whenever $A^{\ast}A=P$ and $AA^{\ast}=Q$.
\end{lemma}

Proof. Let $f:\mathbb{R\rightarrow R}$ be continuous defined by $f(t)=\left\{
\begin{array}
[c]{ll}%
0, & 0\leq t\leq\frac{1}{4}\\
\frac{1}{\sqrt{t}}, & \frac{3}{4}\leq t\leq\frac{5}{4}%
\end{array}
\right.  ,$ and define $\psi(P,Q,A)=f(QAPA^{\ast}Q)QAP$. \hfill$\Box$

\vspace{0.2cm}

The following is our main theorem in this section. Suppose $\mathcal{A}$ is a
unital C*-algebra generated by partial isometries $V_{1},\ldots,V_{n}$ and
C*$\left(  V_{1}^{\ast}V_{1},\ldots,V_{n}^{\ast}V_{n}\right)  $ and C*$\left(
V_{1}V_{1}^{\ast},\ldots,V_{n}V_{n}^{\ast}\right)  $ are both (weakly)
semiprojective or C*$\left(  V_{1}^{\ast}V_{1},\ldots,V_{n}^{\ast}V_{n}%
,V_{1}V_{1}^{\ast},\ldots,V_{n}V_{n}^{\ast}\right)  $ is (weakly)
semiprojective. Does it follow that $\mathcal{A}$ is weakly semiprojective? We
prove this is true when the only relations on $V_{1},\ldots,V_{n}$ are those
on $V_{1}^{\ast}V_{1},\ldots,V_{n}^{\ast}V_{n},V_{1}V_{1}^{\ast},\ldots
,V_{n}V_{n}^{\ast}$.

\begin{theorem}
\label{theorem,glue together, projections} The following are true:

(1)\ Suppose $C^{\ast}\left(  P_{1},\ldots,P_{n}|\ \varphi\right)  $ and
$C^{\ast}\left(  Q_{1},\ldots,Q_{n}|\psi\right)  $ are (weakly)
semiprojective, where $P_{1},Q_{1},\ldots,P_{n},Q_{n}$ are projections. Then
the universal C*-algebra $\mathcal{A}=C^{\ast}(V_{1},\ldots,V_{n}|\Phi)$ with
the relation $\Phi$ defined by
\begin{align*}
\Phi(V_{1},\ldots,V_{n})  =& \varphi\left(  V_{1}^{\ast}V_{1},\ldots
,V_{n}^{\ast}V_{n}\right)  ^{\ast}\varphi\left(  V_{1}^{\ast}V_{1}%
,\ldots,V_{n}^{\ast}V_{n}\right) +\\
&  +\psi\left(  V_{1}V_{1}^{\ast},\ldots,V_{n}V_{n}^{\ast}\right)  ^{\ast}%
\psi\left(  V_{1}V_{1}^{\ast},\ldots,V_{n}V_{n}^{\ast}\right)
\end{align*}
is (weakly) semiprojective.

(2)\ If $C^{\ast}(P_{1},\ldots,P_{n},Q_{1},\ldots,Q_{n}|\phi)$ is
(weakly) semiprojective, where
$P_{1},\ldots,P_{n},Q_{1},\ldots,Q_{n}$ are projections, then
$C^{\ast}(V_{1},\ldots,V_{n}|\Psi)$ is (weakly) semiprojective,
where the relation $\Psi$ is defined by
\begin{align*}
& \Psi(V_{1},\ldots,V_{n})\\
 =& \phi\left(  V_{1}^{\ast}V_{1},\ldots,V_{n}^{\ast}V_{n},V_{1}V_{1}^{\ast
},\ldots,V_{n}V_{n}^{\ast}\right) ^{\ast}\phi\left(  V_{1}^{\ast}V_{1}%
,\ldots,V_{n}^{\ast}V_{n},V_{1}V_{1}^{\ast},\ldots,V_{n}V_{n}^{\ast}\right)
\end{align*}

\end{theorem}

Proof. \ (1)\ Let $\varphi_{1}, \ldots, \varphi_{n}$ be weakly semiprojective
approximating functions for $\varphi$, and $\psi_{1}, \ldots, \psi_{n}$ be
weakly semiprojective approximating functions for $\psi$.

Define the functions $\Phi_{1}, \ldots, \Phi_{n}$ by
$$\Phi_{i}(V_{1}, \ldots, V_{n})=\phi(\varphi_{i}(P_{1}, \ldots,
P_{n}), \psi_{i}(Q_{1}, \ldots, Q_{n}), V_{i}),$$ where $\phi$ is
the noncommutative continuous function in Lemma \ref{lemma,for
theorem in section3}.

Given any $\varepsilon>0$. Let $\delta_{0}$ be defined in Lemma
\ref{lemma,for theorem in section3} corresponding to $\varepsilon$.

Since $C^{\ast}\left(  P_{1},\ldots,P_{n}|\ \varphi\right)  $ is weakly
semiprojective, there exists $\delta_{1}>0$, such that for any operators
$A_{1}, \ldots, A_{n}$ with $\|\varphi(A_{1}^{*}A_{1}, \ldots, A_{n}^{*}%
A_{n})\|<\delta_{1}$, we have, for $1\leq i\leq n$,

(a)\ $\varphi_{i}(A_{1}%
^{*}A_{1}, \ldots, A_{n}^{*}A_{n})$ is projection and

(b)\ $\|A_{i}^{*}%
A_{i}-\varphi_{i}(A_{1}^{*}A_{1}, \ldots, A_{n}^{*}A_{n})\|<\delta_{0}$.
\newline
since $C^*(Q_1, \ldots, Q_n|\psi)$ is weakly semiprojective, there
exists $\delta_{2}>0$, such that if $\|\varphi(A_{1}A_{1}^{*},
\ldots,
A_{n}A_{n}^{*})\|<\delta_{2}$, then, for $1\leq i\leq n$,

(a$'$)\ $\psi_{i}%
(A_{1}A_{1}^{*}, \ldots, A_{n}A_{n}^{*})$ is projection and

(b$'$)\ $\|A_{i}A_{i}%
^{*}-\psi_{i}(A_{1}A_{1}^{*}, \ldots, A_{n}A_{n}^{*})\|<\delta_{0}$.

Put $\varphi_{i}(A_{1}^{*}A_{1}, \ldots, A_{n}^{*}A_{n}), \psi_{i}(A_{1}%
A_{1}^{*}, \ldots, A_{n}A_{n}^{*}), A_{i}$ to $P, Q, A$ in Lemma
\ref{lemma,for theorem in section3}, we have that
\[
\phi(\varphi_{i}(A_{1}^{*}A_{1}, \ldots, A_{n}^{*}A_{n}),\psi_{i}(A_{1}%
A_{1}^{*}, \ldots, A_{n}A_{n}^{*}), A_{i}) \ \left(=\Phi_{i}(A_{1},
\ldots, A_{n})\right)
\]
is a partial isometry from $\varphi_{i}(A_{1}^{*}A_{1}, \ldots, A_{n}^{*}%
A_{n})$ to $\psi_{i}(A_{1}A_{1}^{*}, \ldots, A_{n}A_{n}^{*}).$

Let $\delta=\min\{\delta_{1}^{2}, \delta_{2}^{2}\}$. We prove that $\Phi_{1},
\ldots, \Phi_{n}$ are weakly semiprojective approximating functions for
$\Phi.$

Use the similar idea and Lemma \ref{lemma,for theorem in section3}, we can
prove the weakly semiprojective case.

(2)\ Similar to the proof of Part (1).\hfill$\Box$

\begin{example}
Suppose
\begin{align*}
{\cal M}_{2}(C)  & =C^{\ast}\left(  P_{1}=\left(
\begin{array}
[c]{ll}%
1 & 0\\
0 & 0
\end{array}
\right)  ,P_{2}=\left(
\begin{array}
[c]{ll}%
\frac{1}{2} & \frac{1}{2}\\
\frac{1}{2} & \frac{1}{2}%
\end{array}
\right)  ,P_{3}=\left(
\begin{array}
[c]{ll}%
1 & 0\\
0 & 1
\end{array}
\right)  \right) \\
& = C^{\ast}\left(  P_{1},P_{2},P_{3}|\varphi\right)
\end{align*}
and
\begin{align*}
{\cal M}_{3}(C)  & =C^{\ast}\left(  Q_{1}=\left(
\begin{array}
[c]{lll}%
1 & 0 & 0\\
0 & 0 & 0\\
0 & 0 & 0
\end{array}
\right)  ,Q_{2}=\left(
\begin{array}
[c]{lll}%
0 & 0 & 0\\
0 & 1 & 0\\
0 & 0 & 0
\end{array}
\right)  ,Q_{3}=\left(
\begin{array}
[c]{lll}%
\frac{1}{3} & \frac{1}{3} & \frac{1}{3}\\
\frac{1}{3} & \frac{1}{3} & \frac{1}{3}\\
\frac{1}{3} & \frac{1}{3} & \frac{1}{3}%
\end{array}
\right)  \right) \\
& =  C^{\ast}\left(  Q_{1},Q_{2},Q_{3}|\psi\right)  .
\end{align*}
Then the universal C*-algebra generated by partial isometries $V_{1}%
,V_{2},V_{3}$ such that $$\varphi(
V_{1}^{\ast}V_{1},V_{2}^{\ast}V_{2} ,V_{3}^{\ast}V_{3}) =0$$ and
$$\psi\left( V_{1}V_{1}^{\ast},V_{2}V_{2}^{\ast
},V_{3}V_{3}^{\ast}\right)  =0$$ is semiprojective.
\end{example}

We also can apply our results to the generalized version of the
noncommutative unitary construction of K. McClanahan \cite{M}.

\begin{proposition}
If $\mathcal{A}=C^{\ast}\left(  x_{1},\ldots x_{n}|\varphi\right)  $
is (weakly) semiprojective, then the universal C*-algebra
$\mathcal{E}$ generated by $\left\{  a_{ijk}:1\leq i,j\leq m,1\leq
k\leq n\right\}$, subject to $\varphi\left(  \left( a_{ij1}\right)
,\ldots,\left(  a_{ijn}\right) \right)  =0$ is (weakly)
semiprojective.
\end{proposition}

Proof. Suppose $\varphi_{1}, \ldots, \varphi_{n}$ are (weakly) semiprojective
approximating functions for $\varphi$. Define functions $\{\Phi_{i,j,k}: 1\leq
i,j\leq m,1\leq k\leq n\}$ by
\[
\Phi_{i,j,k}\left(  \{a_{s,t,l}\}_{s,t,l}\right)  =f_{i,j}\left(  \varphi
_{k}\left(  \left(  a_{s,t,1}\right)  ,\ldots,\left(  a_{s,t,n}\right)
\right)  \right)  ,
\]
where $f_{i,j}: {\mathcal{M}}_{m}(\mathbb{C)\mapsto C}$ such that for any
$m\times m$ matrix $A$, $A=\left(  f_{i,j}(A)\right)  $. It is clear that
$\{\Phi_{i,j,k}: 1\leq i,j\leq m,1\leq k\leq n\}$ are (weakly) semiprojective
approximating functions for $\Phi.$\hfill$\Box$

\begin{corollary}
Suppose $C^{\ast}\left(  P_{1},\ldots,P_{n}|\varphi\right)  $ and $C^{\ast
}\left(  Q_{1},\ldots,Q_{n}|\psi\right)  $ are (weakly) semiprojective, where
$P_{1},Q_{1},\ldots,P_{n},Q_{n}$ are projections. Suppose $m$ is a positive
integer and $\mathcal{A}$ is the universal C*-algebra generated by
$\{a_{ijk}:$$1\leq i,j\leq m,1\leq k\leq n\}$ subject to
\[
\varphi\left(  \left(  a_{ij1}\right)  ^{\ast}\left(  a_{ij1}\right)
,\ldots,\left(  a_{ijn}\right)  ^{\ast}\left(  a_{ijn}\right)  \right)  =0,
\]
\[
\psi((a_{ij1})(a_{ij1})^{\ast},\ldots,(a_{ijn})(a_{ijn})^{\ast})=0.
\]
Then $\mathcal{A}$ is (weakly) semiprojective.
\end{corollary}

We can also define projectivity in terms of noncommutative continuous
functions. A unital C$^{\ast}$-algebra C*$\left(  b_{1},\ldots,b_{n}%
|\varphi\right)  $ is projective in the unital category if, for any unital
C*-algebra $\mathcal{A}$ and any ideal $\mathcal{J}$ in $\mathcal{A}$, and any
$x_{1},\ldots,x_{n}\in{\mathcal{A}}/{\mathcal{J}}$ with $\varphi(x_{1}%
,\ldots,x_{n})=0$, there exist elements $a_{1},\ldots,a_{n}$ in
$\mathcal{A}$, such that $x_{i}=a_{i}+{\mathcal{J}}$ and
$\varphi(a_{1},\ldots,a_{n})=0$.

\begin{remark}
(1)\  The universal C*-algebra generated by $A$ such that
$\left\Vert A\right\Vert \leq r$ is projective. The universal
C*-algebra generated by $\left\{  A_{n}\right\}_{n=1}^{\infty}$ such
that $\left\Vert A_{n}\right\Vert \leq r_{n}$ for some numbers
$r_n$, is projective. Thus every separable unital C*-algebra is
isomorphic to $\mathcal{A}/\mathcal{J}$ , where $\mathcal{A}$ is a
projective C$^*$-algebra and $\cal J$ is an ideal of $\cal A$.

(2)\ If $\{{\mathcal{A}}_{n}\}_{n=1}^{\infty}$ is a sequence of
projective C$^*$-algebras, then the free product
$\ast_{n}{\mathcal{A}}_{n}$ is projective.

(3)\ If $\left\{  \mathcal{A}_{n}\right\}_{n=1}^{\infty}$ is a
sequence of separable unital semiprojective algebras, then
$\ast_{n}\mathcal{A}_{n}$ may not be weakly semiprojective. For
example,
${\mathcal{M}}_{2}(\mathbb{C})\ast{\mathcal{M}}_{3}(\mathbb{C)\ast
{\mathcal{M}}}_{4}\mathbb{(C)\ast\cdots}$ is not weakly
semiprojective, but each ${\mathcal{M}}_{n}(\mathbb{C})$ is
semiprojective.
\end{remark}

Although weakly semiprojective C*-algebras need not be finitely generated,
identity representation on such algebras must be a pointwise limit of
representations into finitely generated subalgebras.

\begin{proposition}
\label{DL}Suppose $\mathcal{A}$ is separable and weakly
semiprojective and $\left\{\mathcal{A}_{n}\right\}_{n=1}^{\infty}$
is an increasing sequence of finitely generated C*-subalgebras whose
union is dense in $\mathcal{A}$. Then there is
a positive integer $N$ and unital *-homomorphisms $\pi_{n}%
:\mathcal{A\rightarrow A}_{n}$ for all $n\geq N$ such that%
\[
\left\Vert x-\pi_{n}\left(  x\right)  \right\Vert \rightarrow0
\]
for every $x\in\mathcal{A}$.
\end{proposition}

Proof. It follows from the hypothesis that dist$\left(  x,\mathcal{A}%
_{n}\right)  \rightarrow0$ for every $x\in\mathcal{A}$. Thus, for
each $x\in\mathcal{A}$, there is a $x_n\in\mathcal{A}_{n}$ such that
$\left\Vert x-x_n\right\Vert \leq$dist$\left(
x,\mathcal{A}_{n}\right) +\frac{1}{n}$. Define a unital
*-homomorphism $\pi:\mathcal{A}\rightarrow\prod_{1}^{\infty}
\mathcal{A}_{n}/\oplus _{1}^{\infty} \mathcal{A}_{n}$ by
\[
\pi\left(  x\right)  =\left[  \left\{x_n\right\} \right] .
\]
The desired result follows easily from the weak semiprojectivity of
$\mathcal{A}$.\hfill$\Box$

\begin{definition}
A unital $C^{\ast}$-algebra $\mathcal{A}$ is called GCR if for any
irreducible representation $\pi$ from $\mathcal{A}$ to
$\mathcal{B}(\mathcal{H})$,
$\mathcal{K}(\mathcal{H})\subseteq\pi(\mathcal{A})$.
\end{definition}

\begin{lemma}
If $\mathcal{A}$ is a unital GCR C$^{\ast}$-algebra, then there exists a
positive integer $n$ and a representation $\pi:{\mathcal{A}}\rightarrow
{\mathcal{M}}_{n}(\mathbb{C)}$ that is onto.
\end{lemma}

Proof. Suppose $\mathcal{J}$ is a maximal ideal in $\mathcal{A}$. Then
${\mathcal{A}}/{\mathcal{J}}$ is a simple C$^{*}$-algebra. Let $\pi:
{\mathcal{A}}/{\mathcal{J}}\rightarrow{\mathcal{B}(\mathcal{H})}$ be an
irreducible representation. Then $\pi\left(  {\mathcal{A}}/{\mathcal{J}%
}\right)  ^{\prime}=\mathbb{C}1$. It follows that ${\mathcal{K}(\mathcal{H}%
)}\subseteq\pi\left(  {\mathcal{A}}/{\mathcal{J}}\right)  $ is a closed ideal,
therefore $\mathcal{H}$ is finite-dimensional. \hfill$\Box$

\vspace{0.1cm}

From the above lemma, it is not hard to see that if $\mathcal{A}$ is a simple
infinite-dimensional C$^{\ast}$-algebra, then $\mathcal{A}$ cannot be a
subalgebra of a GCR algebra.

\begin{corollary}
If $\mathcal{A}$ is a unital simple infinite-dimensional C*-algebra that is a
subalgebra of a direct limit of subalgebras of GCR C*-algebras, then
$\mathcal{A}$ is not weakly\ semiprojective.
\end{corollary}

Proof. It follows from the proof of Proposition \ref{DL} that there is a
*-homomorphism $\pi:\mathcal{A}\rightarrow\prod_{1}^{\infty}
\mathcal{A}_{n}/\oplus_{1}^{\infty} \mathcal{A}_{n}$, where $\{{\cal
A}_n\}_{n=1}^{\infty}$ as defined in Proposition \ref{DL}. Assume
via contradiction that $\mathcal{A}$ is weakly semiprojective. Then
there is a representation $\pi_{n}:\mathcal{A\rightarrow A}_{n}$ for
some positive integer $n$. Since $\mathcal{A}_{n}$ is a subalgebra
of a GCR algebra, it follows that $\mathcal{A}_{n},$ and hence
$\mathcal{A}$, has a finite-dimensional representation. Since
$\mathcal{A}$ is simple, every representation of $\mathcal{A}$ is
one-to-one, which implies that $\mathcal{A}$ is finite-dimensional,
a contradiction.\hfill$\Box$

\begin{remark}
The Cuntz algebra is weakly semiprojective, hence, the Cuntz algebra cannot be
embedded into a direct limit of subalgebras of GCR C*-algebras. The irrational
rotation algebra ${\mathcal{A}}_{\theta}$ is not weakly semiprojective, since
it can be embedded into the direct limit of subalgebras of GCR C*-algebras.
\end{remark}

We conclude this section with an observation concerning the reduced
free group C$^{*}$-algebra, $C_{r}^{\ast}\left(
\mathbb{F}_{n}\right)$.

\begin{proposition}
$C_{r}^{\ast}\left(  \mathbb{F}_{n}\right)  $ is not weakly semiprojective.
\end{proposition}

Proof. U. Haagerup and S. Thorbj\o rnsen \cite{Haag} proved that
there is a map $$\pi:C_{r}^{\ast}\left(  \mathbb{F}_{n}\right)
\rightarrow\prod_{n\geq
1}\mathcal{M}_{n}\left(  \mathbb{C}\right)  /\oplus_{n\geq1}\mathcal{M}%
_{n}\left(  \mathbb{C}\right)  .$$ Since $C_{r}^{\ast}\left(  \mathbb{F}%
_{n}\right)  $ is infinite-dimensional and simple,
$C_{r}^{\ast}\left( \mathbb{F}_{n}\right)$ has no finite-dimensional
representation. Hence $C_{r}^{\ast}\left( \mathbb{F}_{n}\right)  $
cannot be weakly semiprojective.\hfill$\Box$

\section{Finite Von Neumann Algebras and trace norms}

When we talked about weak semiprojectivity in C*-algebras we described it in
terms of mappings into algebras $\prod_{1}^{\infty}{\mathcal{B}}_{n}%
/\oplus_{1}^{\infty}{\mathcal{B}}_{n}$ being ``liftable". There is another way
to describe this by replacing the $\prod_{1}^{\infty}{\mathcal{B}}_{n}%
/\oplus_{1}^{\infty}{\mathcal{B}}_{n}$ construction with ultraproducts.

Suppose $\mathbb{I}$ is an infinite set and $\omega$ is an ultrafilter on
$\mathbb{I}$, i.e., $\omega$ is a family of subset of $\mathbb{I}$ such that
\newline(1) $\emptyset\notin\omega$ \newline(2) If $A,B\in\omega$, then $A\cap
B\in\omega$ \newline(3) For every subset $A$ in $\mathbb{I}$, either
$A\in\omega$ or $\mathbb{I}\setminus A\in\omega$.

One example of an ultrafilter is obtained by choosing an element
$\iota$ in $\mathbb{I}$ and letting $\omega$ be the collection of
all subsets of $\mathbb{I}$ that contain $\iota$. Such an
ultrafilter is called \emph{principle} ultrafilter and ultrafilter
not of this form are called \emph{free}. We call an ultrafilter
$\omega$ \emph{nontrivial} if it is free and there is a sequence
$\left\{  E_{n}\right\}_{n=1}^{\infty}$ in $\omega$ whose
intersection is empty. We can always choose $E_{1}=\mathbb{I}$ and,
by replacing $E_{n}$ with $\cap_{k=1}^{n}E_{k}$, we can assume that
$\left\{ E_{n}\right\}_{n=1}^{\infty}$ is decreasing. Throughout
this paper we will only use nontrivial ultrafilters.

Suppose $\{{\mathcal{A}}_{i}:i\in\mathbb{I\}}$ is a family of C*-algebras and
$\omega$ is a nontrivial ultrafilter on $\mathbb{I}$. Then
\[
{\mathcal{J}}=\left\{  \{A_{i}\}\in\prod_{i\in\mathbb{I}}{\mathcal{A}}%
_{i}:\ \lim_{i\rightarrow\omega}\left\Vert A_{i}\right\Vert =0\right\}
\]
is a norm-closed two-sided ideal in
$\prod_{i\in\mathbb{I}}{\mathcal{A}}_{i},$ and we call the quotient
the \emph{C*-algebraic ultraproduct} of the $\mathcal{A}_{i}$'s and
denote it by $\prod^{\omega}{\mathcal{A}}_{i}$. For an introduction
to ultraproducts see \cite{Don1}. It is easily verified that a
C*-algebra $\mathcal{A}$ is weakly semiprojective if and only if,
given a unital *-homomorphism $\pi:\mathcal{A}\rightarrow$
$\prod^{\omega }{\mathcal{A}}_{i}$ there are functions
$\pi_{i}:\mathcal{A}\rightarrow \mathcal{A}_{i}$ for each
$i\in\mathbb{I}$ such that, eventually along $\omega,$ $\pi_{i}$ is
a unital *-homomorphism and such that, for ever $a\in\mathcal{A}$,
\[
\pi\left(  a\right)  =[\left\{  \pi_{i}\left(  a\right)  \right\} ] _{\omega}.
\]

We now want to look at analogue of weak semiprojectivity for finite
von Neumann algebras with faithful tracial states. Suppose
$\mathcal{A}$ is a C*-algebra with a tracial state $\tau$. As in the
GNS construction there is a seminorm $\left\Vert \cdot\right\Vert
_{2, \tau}$ on $\mathcal{A}$ defined by $\left\Vert a\right\Vert
_{2,\tau}=\tau\left( a^{\ast}a\right)  ^{1/2}$. More generally, if
$1\leq p<\infty$, we define $\Vert a\Vert_{p,\tau}=(\tau((a^{\ast
}a)^{p/2}))^{1/p}$. Since $C^{\ast}(a^{\ast}a)$ is isomorphic to
$C(X)$, where $X$ is the spectrum of $a^{\ast}a$, there is a
probability measure $\mu$ such that
$\tau(f(a^{\ast}a))=\int_{X}fd\mu$ for every $f\in C(X)$. Thus, for
$1\leq p<\infty$, \begin{equation}\label{equation,1}\Vert
a\Vert_{p,\tau}=0\ \  \mbox{if and only if}\ \  \Vert a\Vert_{2,
\tau}=0.\end{equation} If there is no confusion, we can simply use
$\|\cdot\|_2$ and $\|\cdot\|_p$ to denote $\|\cdot\|_{2,\tau}$ and
$\|\cdot\|_{p,\tau}$ respectively.

Suppose $\left\{  \left(  \mathcal{A}_{i},\tau_{i}\right)  :i\in
\mathbb{I}\right\}  $ is a family of C*-algebras $\mathcal{A}_{i}$ with
tracial states $\tau_{i}$. We can define a trace $\rho$ on $\prod
_{i\in\mathbb{I}}{\mathcal{A}}_{i}$ by
\[
\rho\left(  \left\{  a_{i}\right\}  \right)  =\lim_{i\rightarrow\omega}%
\tau_{i}\left(  a_{i}\right)  .
\]

The set $\mathcal{J}_{2}=\left\{  \{A_{i}\}\in\prod_{i\in\mathbb{I}%
}{\mathcal{A}}_{i}:\ \lim_{i\rightarrow\omega}\left\Vert A_{i}\right\Vert
_{2}=0\right\}  $ is a closed two-sided ideal in $\prod_{i\in\mathbb{I}%
}{\mathcal{A}}_{i}$, and we call the quotient $\left(  \prod_{i\in\mathbb{I}%
}{\mathcal{A}}_{i}\right)  /\mathcal{J}_{2}$ the \emph{tracial ultraproduct}
of the $\mathcal{A}_{i}$'s, and we denote it by $\prod^{\omega}\left(
\mathcal{A}_{i},\tau_{i}\right)  $. There is a natural faithful trace $\tau$
on $\prod^{\omega}\left(  \mathcal{A}_{i},\tau_{i}\right)  $ defined by
\[
\tau\left(  \lbrack\left\{  a_{i}\right\}  ]_{\omega}\right)  =\lim
_{i\rightarrow\omega}\tau_{i}\left(  a_{i}\right)  .
\]
$\prod^{\omega}\left(  \mathcal{A}_{i},\tau_{i}\right)  $ is the
representation of $\prod_{i\in\mathbb{I}}{\mathcal{A}}_{i}$ using
the GNS construction with $\rho$. By Equation (\ref{equation,1}) we
see that $\prod^{\omega}\left( \mathcal{A}_{i},\tau_{i}\right)
=\left( \prod_{i\in\mathbb{I}}{\mathcal{A}}_{i}\right)
/\mathcal{J}_{p}$, where
$\mathcal{J}_{p}=\left\{  \{A_{i}\}\in\prod_{i\in\mathbb{I}}{\mathcal{A}}%
_{i}:\ \lim_{i\rightarrow\omega}\left\Vert A_{i}\right\Vert _{p}=0\right\}  .$
One immediate consequence of $\mathcal{J}_{p}=\mathcal{J}_{2}$ is the fact
that, on a bounded subset of any $\left(  \mathcal{A},\tau\right)  $ the norms
$\left\Vert \cdot\right\Vert _{2}$ and $\left\Vert \cdot\right\Vert _{p}$
generate the same topology.

It was shown by S. Sakai \cite{sakai} that a tracial ultraproduct of
finite factors is a von Neumann algebra and is, in fact, a factor.
However, it is true that any tracial ultraproduct of von Neumann
algebras is a von Neumann algebra. Here we prove that any tracial
ultraproduct of C*-algebras is a von Neumann algebra.

\begin{theorem}
\label{theorem, saikai} Suppose
$\{\mathcal{A}_{i}\}_{i\in\mathbb{I}}$ is a family of
C$^{*}$-algebras with a tracial state $\tau_{i}$ on each
$\mathcal{A}_{i} $ and $\omega$ is a nontrivial ultrafilter on
$\mathbb{I}$. Then the tracial ultraproduct $\prod^{\omega}\left(
{\mathcal{A}}_{i},\tau _{i}\right)  $ of
$\{\mathcal{A}_{i}\}_{i\in\mathbb{I}}$ is a von Neumann algebra.
\end{theorem}

Proof. Let $\mathcal{A}=\prod^{\omega}\left(  {\mathcal{A}}_{i},\tau
_{i}\right)  $. Note that $Ball(\overline{\mathcal{A}}^{\ast-SOT}%
)=\overline{Ball(\mathcal{A})}^{\Vert\cdot\Vert_{2}},$ i.e., the unit ball of
$\overline{\mathcal{A}}^{\ast-SOT}$ is equal to the $\Vert\cdot\Vert_{2}$
closure of the unite ball of $\mathcal{A}$.

Suppose $T\in\overline{Ball(\mathcal{A})}^{\|\cdot\|_{2}}$. Then for
any positive integer $n$, there exists $A_{n}\in Ball(\mathcal{A})$
such that $\|T-A_{n}\|_{2}\leq\frac{1}{4^{n}}.$ Write
$A_{n}=[\{A_{ni}\}]_{\omega}$ with each $A_{ni}\in
Ball(\mathcal{A}_{i})$.

Since $\omega$ is nontrivial, there is a family $\{E_{n}\}$ of elements of
$\omega$ such that
\[
I=E_{1}\supseteq E_{2}\supseteq\cdots\ \ \mbox{and}\ \ \cap_{n}E_{n}%
=\emptyset.
\]
Let
\[
F_{n}=\{i\in E_{n}: \forall1\leq k\leq n, \|A_{ki}-A_{ni}\|_{2}<\frac{1}%
{4^{n}}+\frac{1}{4^{k}}\}.
\]
Let $X_{i}=A_{ki}$ for $i\in F_{k}/F_{k+1}$. For any $i\in F_{n}$, there
exists some $k\geq n$ such that $i\in F_{k}/F_{k+1}$ and
\[
\|A_{ni}-X_{i}\|_{2}=\|A_{ni}-A_{ki}\|_{2}\leq\frac{1}{4^{n}}+\frac{1}{4^{k}%
}\leq\frac{2}{4^{n}}\leq\frac{1}{2^{n}}.
\]
Let $X=[\{X_{i}\}]_{\omega}$. Then $X\in Ball (\mathcal{A})$ and
\[
\|A_{n}-X\|_{2}\leq\frac{1}{2^{n}}.
\]
Hence $T=X\in Ball (\mathcal{A}).$ This implies that ${\cal
A}=\prod^{\omega}\left( {\mathcal{A}}_{i},\tau _{i}\right) $ is a
von Neumann algebra.\hfill$\Box$

\vspace{0.2cm}

The next Theorem gives a generalization of Lin's theorem for
$\left\Vert \cdot\right\Vert _{p}$ on C*-algebras with trace. When
$p=2$, it was proved for finite factors in \cite{Don2}.

\begin{theorem}
\label{corollary,shanghai}For every $\varepsilon>0$ and every $1\leq
p<\infty,$ there exists $\delta>0$ such that, for any
C$^{\ast}$-algebra $\mathcal{A}$ with trace $\tau$, and
$A_{1},\ldots,A_{n}\in$ball$\left( \mathcal{A}\right)  $ with $\Vert
A_{j}A_{j}^{\ast}-A_{j}^{\ast}A_{j}\Vert _{p}<\delta$ and $\Vert
A_{j}A_{k}-A_{k}A_{j}\Vert_{p}<\delta$, there exists
$B_{1},\ldots,B_{n}\in$ball$\left(  \mathcal{A}\right)  $ so that
$B_{j}B_{j}^{\ast}=B_{j}^{\ast}B_{j}$, $B_{j}B_{k}=B_{k}B_{j}$ and
$\sum _{j=1}^{n}\Vert A_{j}-B_{j}\Vert_{p}<\varepsilon$.
\end{theorem}

Proof. Assume the statement is false. Then there is an $\varepsilon>0$ such
that, for every positive integer $k$, there is a unital C*-algebra
$\mathcal{A}_{k}$ with trace $\tau_{k}$ and elements $A_{k,1},\ldots,A_{k,n}$
with $\Vert A_{k,j}A_{k,j}^{\ast}-A_{k,j}^{\ast}A_{k,j}\Vert_{p}<\frac{1}{k}$
and $\Vert A_{k,j}A_{k,i}-A_{k,i}A_{k,j}\Vert_{p}<\frac{1}{k}$, so that for
all $B_{1},\ldots,B_{n}\in$ball$\left(  \mathcal{A}\right)  $ with $B_{j}%
B_{j}^{\ast}=B_{j}^{\ast}B_{j}$ and $B_{j}B_{k}=B_{k}B_{j}$ we have
$\sum_{j=1}^{n}\Vert A_{j}-B_{j}\Vert_{p}^{p}\geq\varepsilon.$ The
tracial ultraproduct
$\mathcal{A=}\prod^{\omega}\mathcal{A}_{i}=\left(  \prod
_{i\in\mathbb{I}}{\mathcal{A}}_{i}\right)  /\mathcal{J}_{p}$ is a
von Neumann algebra and $\{A_{j}=\left[  \left\{  A_{k,j}\right\}
\right]  _{\omega}: 1\leq j\leq n\}$ is a family of commuting normal
operators. Hence, by the proof of Theorem 5.5 in \cite{BDF}, there
is a selfadjoint operator $C\in\mathcal{A}$ and bounded continuous
functions $f_{1},\ldots ,f_{n}:\mathbb{R}\rightarrow\mathbb{C}$ such
that $A_{j}=f_{j}\left( C\right)  $ for $1\leq j\leq n.$ Write
$C=\left[  \left\{ C_{k}\right\}
\right]  _{\omega}$ with each $C_{k}=C_{k}^{\ast}$. Define $B_{k,j}%
=f_{j}\left(  C_{k}\right)  $ for $1\leq j\leq n$ and
$k\in\mathbb{N}$. Then $A_{j}=\left[  \left\{  B_{k,j}\right\}
\right]  _{\omega}$ for $1\leq j\leq n$ and $\left\{  B_{k,j}:1\leq
j\leq n\right\}  $ is a family of commuting normal operators. So
\[
\varepsilon\leq\lim_{k\rightarrow\omega}\sum_{j=1}^{n}\Vert A_{k,j}%
-B_{k,j}\Vert_{p}=\sum_{j=1}^{n}\Vert A_{j}-f_{j}\left(  C\right)  \Vert
_{p}=0,
\]
which is a contradiction.\hfill$\Box$

\begin{remark}
Suppose $K$ is a compact nonempty subset of $\mathbb{C}$ that is a
continuous image of $[0,1]$. It follows from Proposition 39 in
\cite{Don} that there is a noncommutative continuous function
$\alpha$ such that, for every operator $T$ with $\Vert T\Vert\leq1$
we have $\alpha(T)$=0 if and only if $T$ is normal and the spectrum
of $T$ is contained in $K$. If, in Corollary
\ref{corollary,shanghai}, we add the condition that $\Vert
\alpha(A_{1})\Vert_{2}<\delta$, then we can choose $B_{1}$ so that
its spectrum is contained in $K$. In particular, if we add $\Vert
1-A_{1}^{\ast}A_{1}\Vert_{2}<\delta$, we can choose $B_{1}$ to be
unitary.
\end{remark}

The next theorem shows that, unlike in the C*-algebra case, commutative
C*-algebras are ``weakly semiprojective" in the ``diffuse von Neumann algebra" sense.

\begin{theorem}
\label{theorem,diffuse}Suppose, for $i\in \Bbb I$, $\mathcal{M}_{i}$
is a diffuse von Neumann algebra with faithful trace $\tau_i$ and
$\mathcal{A}$ is a commutative countably generated von Neumann
subalgebra of the ultraproduct $\prod^{\omega}\left(
\mathcal{M}_{i},\tau_{i}\right) .$ Then, for every $i,$ there is a
trace-preserving *-homomorphism $\pi_{i}:\mathcal{A}\rightarrow
\mathcal{M}_{i}$ such that, for every $a\in\mathcal{A},$
\[
a=[\left\{  \pi_{i}\left(  a\right)  \right\} ] _{\omega}.
\]

\end{theorem}

Proof. Suppose $P$ is a projection in $\prod^{\omega}\left(  \mathcal{M}%
_{i},\tau_{i}\right)  $. It is well-known that $P$ can be written as
$P=[\{A_{i}\}]_{\omega}$ with each $A_{i}$ a projection. Since $\tau
(A_{i})\rightarrow\tau(P)$ as $i\rightarrow\omega$ and since each
$\mathcal{M}_{i}$ is diffuse, we can, for each $i$, find a
projection $P_{i}\in\mathcal{M}_{i}$ so that $\tau_{i}\left(
P_{i}\right) =\tau\left(  P\right)$ and either $P_{i}\leq A_{i}$ or
$A_{i}\leq P_{i}$. Since $\left\Vert A_{i}-P_{i}\right\Vert
_{2}=\sqrt{\left\vert \tau_{i}\left( P_{i}\right)  -\tau_{i}\left(
A_{i}\right)  \right\vert }\rightarrow0$, we have $P=[\left\{
P_{i}\right\} ]_{\omega}$. Hence, every projection in
$\prod^{\omega}\left( \mathcal{M}_{i},\tau_{i}\right)  $ can be
lifted to projections with the same trace.

Next suppose $P=[\left\{  P_{i}\right\} ]_{\omega},Q=[\left\{  Q_{i}\right\}
]_{\omega}$ are projections in $\prod^{\omega}\left(  \mathcal{M}_{i},\tau
_{i}\right)  $ such that $P\leq Q,$ and, for every $i$, $P_{i}\leq Q_{i}$ and
$\tau_{i}\left(  P_{i}\right)  =\tau\left(  P\right)  $ and $\tau_{i}\left(
Q_{i}\right)  =\tau\left(  Q\right)  $. Suppose $E$ is a projection in
$\prod^{\omega}\left(  \mathcal{M}_{i},\tau_{i}\right)  $ and $P<E<Q$.
Applying what we just proved to the projection $E-P$ in the ultraproduct

\begin{center}
$\left(  Q-P\right)  \left(  \prod^{\omega}\left(  \mathcal{M}_{i},\tau
_{i}\right)  \right)  \left(  Q-P\right)  =\prod^{\omega}\left(  Q_{i}%
-P_{i}\right)  \mathcal{M}_{i}\left(  Q_{i}-P_{i}\right)  ,$
\end{center}

we can find projections $E_{i}\in\mathcal{M}_{i}$ so that $P_{i}\leq
E_{i}\leq Q_{i}$, $\tau_{i}\left(  E_{i}\right)  =\tau\left(
E\right)  $ and $E=[\left\{  E_{i}\right\} ]_{\omega}$. Since
$\mathcal{A}$ is countably generated and commutative, we know from
von Neumann's Theorem that $\mathcal{A}$ is generated by a single
selfadjoint $T$ with $0\leq T\leq1.$ Since $\prod^{\omega}\left(
\mathcal{M}_{i},\tau_{i}\right)  $ is diffuse, the chain $\left\{
\chi_{\lbrack0,s)}\left(  T\right)  :0\leq s\leq1\right\} $ can be
extended to a chain $\left\{  P\left(  t\right)  :t\in\left[
0,1\right]  \right\}  $ such that $\tau\left(  P\left(  t\right)
\right)  =t$ for $0\leq t\leq1$. Repeatedly using the result above
we can find projections $P_i\left(  t\right)$ for each $i$ and each
rational $t\in\left[ 0,1\right]  $ such that $\tau_{i}\left(
P_i\left( t\right)\right)  =t$
and $P\left(  t\right)  =[\left\{  P_i\left(  t\right)\right\} ]_{\omega}%
$, and such that $P_i\left(  s\right)\leq P_i\left(  t\right)$ for
all $i$ and $s\leq t$. Hence, for each $t\in\left[ 0,1\right]  $ and
each $i\in
I,$ we can define%
\[
P_i\left(  t\right)=\sup\left\{  P_i\left(  s\right):s\leq
t,s\in\mathbb{Q}\right\}  =\inf\left\{  P_i\left(  s\right) :s\geq
t,s\in\mathbb{Q}\right\}  .
\]
Then we must have $P\left(  t\right)  =[\left\{  P_i\left(  t\right)
\right\} ] _{\omega}$ for every $t\in\left[  0,1\right] $. For each
$i$, the map $P\left(  t\right)  \mapsto P_i\left(  t\right)$
extends to a trace-preserving *-homomorphism $\rho_{i}:\left\{
P\left(  t\right) :t\in\left[  0,1\right]  \right\}
^{\prime\prime}\rightarrow\mathcal{M}_{i} $, and we can let
$\pi_{i}=\rho_{i}|\mathcal{A}$.\hfill$\Box$

\begin{corollary}
\label{corollary,tracepreserving homomorphism} Suppose $\mathcal{M}_{i}$ is a
diffuse von Neumann algebra for every $i\in\mathbb{I}$ and $\mathcal{A}$
commutative countably generated unital C*-algebra and $\pi:\mathcal{A}%
\rightarrow$ $\prod^{\omega}\left(  \mathcal{M}_{i},\tau_{i}\right)
$ is a unital *-homomorphism. Then, for every $i,$ there is a
*-homomorphism $\pi_{i}:\mathcal{A}\rightarrow\mathcal{M}_{i}$ such
that

\begin{enumerate}
\item $\pi\left(  a\right)  =\left[  \left\{  \pi_{i}\left(  a\right)
\right\}  \right]  _{\omega}$ for every $a\in\mathcal{A}$, and

\item $\tau_{i}\circ\pi_{i}=\tau\circ\pi$ for every $i\in\mathbb{I}$.
\end{enumerate}
\end{corollary}

\begin{corollary}
For every $\varepsilon>0$ there is a positive integer $N$ and
$\delta>0$ such that, for every diffuse finite von Neumann algebra
with trace $\tau$ and every $U\in$ball$\left(  \mathcal{M}\right) ,$
if $\left\vert \tau\left( U^{k}\right)  \right\vert <\delta$ for
$1\leq k\leq N$ with $\left\Vert 1-U^{\ast}U\right\Vert
_{2}<\delta$, then there is a Haar unitary $V\in\mathcal{M}$ such
that $\left\Vert U-V\right\Vert _{2}<\varepsilon$.
\end{corollary}

\begin{remark}
It follows from Theorem \ref{theorem,diffuse} that the hypothesis in Theorem 4
in \cite{Don4} that $w_{1},w_{2},\ldots$ are Haar unitaries can be replaced by
the assumption that they are unitaries. In particular, if $\mathcal{M}$ is a
von Neumann algebra with a faithful trace $\tau$, and $\mathcal{N}$ is a
diffuse subalgebra of $\mathcal{\ M}$, and $\{v_{n}\}$ is a sequence of Haar
unitaries in $\mathcal{M}$ and $\{w_{n}\}$ is a sequence of unitaries in
$\mathcal{N}$, and $\Vert w_{n}-v_{n}\Vert_{2}\rightarrow0$, then there exists
a sequence $\{u_{n}\}$ of Haar unitaries in $\mathcal{\ N}$ such that $\Vert
u_{n}-v_{n}\Vert_{2}\rightarrow0$.
\end{remark}

We next give an analogue of Theorem \ref{theorem,diffuse} with $\mathcal{A}$
hyperfinite instead of commutative, but with each of the $\mathcal{M}_{i}$'s a
II$_{1}$ factor.

\begin{lemma}
\label{lemma,con}\cite{Con} Let $\mathcal{M}$ be a separable factor and
$\omega$ a nontrivial ultrafilter. Let $E=[\{E_{i}\}]_{\omega}$ and
$F=[\{F_{i}\}]_{\omega}$ be equivalent projections in $\prod^{\omega
}\mathcal{M}$ with $E_{i}$'s and $F_{i}$'s projections. Suppose $V$ is a
partial isometry from $E$ to $F$. Then $V=[\{V_{i}\}]_{\omega}$, where $V_{i}$
is a partial isometry from $E_{i}$ to $F_{i}$.
\end{lemma}

\begin{lemma}
\label{lemma,lift trace-preserving homo} Suppose each
$\mathcal{M}_{i}$ is a II$_{1}$ factor with the trace $\tau_i$ and
$\mathcal{A}\subseteq\mathcal{B}$ are finite-dimensional
C*-subalgebras of the ultraproduct $\prod^{\omega}\left(  \mathcal{M}_{i}%
,\tau_{i}\right)  , $ and suppose for every $i,$ there is a trace-preserving
homomorphism $\pi_{i}:\mathcal{A}\rightarrow\mathcal{M}_{i}$ such that, for
every $a\in\mathcal{A},$ $a=[\left\{  \pi_{i}\left(  a\right)  \right\}
]_{\omega}.$ Then for every $i,$ there is a trace-preserving homomorphism
$\rho_{i}:\mathcal{B}\rightarrow\mathcal{M}_{i}$ such that,

(1) for every $b\in\mathcal{B},$
\[
b=[\left\{  \rho_{i}\left(  b\right)  \right\} ]_{\omega},
\]
and

(2) for every $i,$ $\rho_{i}|_{\mathcal{A}}=\pi_{i}$.
\end{lemma}

Proof. To avoid a notational nightmare, we will describe the proof for a
specific example. It will be easy to see how this technique applies
universally. Suppose $\mathcal{A}$ is isomorphic to $\mathcal{M}_{2}%
\oplus\mathcal{M}_{3}$ and $\mathcal{B}$ is isomorphic to $\mathcal{M}%
_{4}\oplus\mathcal{M}_{5}$ where the inclusion $\mathcal{A}\subset\mathcal{B}$
identifies $A\oplus B$ with $\left(  A\oplus A\right)  \oplus\left(  A\oplus
B\right)  $. Let $\left\{  e_{st}:1\leq s,t\leq4\right\}  $ denote matrix
units for $\mathcal{M}_{4}\oplus0$ and $\left\{  f_{st}:1\leq s,t\leq
5\right\}  $ denote matrix units for $0\oplus\mathcal{M}_{5}$. Then
$$\mathcal{S}_{1}=\left\{  e_{11}+e_{33}+f_{11},e_{12}+e_{34}+f_{12}%
,e_{21}+e_{43}+f_{21},e_{22}+e_{44}+f_{22}\right\}  $$ is a set of
matrix units for
$\mathcal{M}_{2}\oplus0$ and $$\mathcal{S}_{2}=\left\{  f_{33},f_{34}%
,f_{35},f_{43},f_{44},f_{45},f_{53},f_{54},f_{55}\right\}  $$ is a
set of matrix units
for $0\oplus\mathcal{M}_{3}$. We have $\pi_{i}$ is defined on $\mathcal{S}%
_{1}\cup\mathcal{S}_{2}.$ We want to extend $\pi_{i}$ to all of the matrix
units for $\mathcal{B}$. However, $\left\{  e_{11}+e_{33}+f_{11},e_{22}%
+e_{44}+f_{22},f_{33},f_{55}\right\}  $ is a commuting family that
is contained in $span\left(  \left\{  e_{ss}:1\leq s\leq4\right\}
\cup\left\{ f_{ss}:1\leq s\leq5\right\}  \right)  ,$ and, using the
techniques in the proof of Theorem \ref{theorem,diffuse}, we can
extend $\pi_{i}$ to a trace-preserving *-homomorphism $\rho_{i}$ on
$span\left(  \left\{ e_{ss}:1\leq s\leq4\right\}  \cup\left\{
f_{ss}:1\leq s\leq5\right\} \right)  $. Using the fact that
\[
e_{11}\left(  e_{12}+e_{34}+f_{12}\right)  =e_{12},
\]
we naturally can define
\[
\rho_{i}\left(  e_{12}\right)  =\rho_{i}\left(  e_{11}\right)  \pi_{i}\left(
e_{12}+e_{34}+f_{12}\right)  .
\]
The definition of $\rho_{i}$ for the remaining matrix units in $\mathcal{B}$
is immediately obtained using Lemma \ref{lemma,con}.\hfill$\Box$

\begin{theorem}
\label{theorem,hyperfinite} If each $\mathcal{M}_{i}$ is a II$_{1}$
factor with the trace $\tau_i$ and $\mathcal{A}$ is a countably
generated hyperfinite von Neumann subalgebra of the ultraproduct
$\prod^{\omega}\left( \mathcal{M}_{i},\tau_{i}\right) ,$ then, for
every $i,$ there is a trace-preserving homomorphism $\pi
_{i}:\mathcal{A}\rightarrow\mathcal{M}_{i}$ such that, for every
$a\in\mathcal{A},$
\[
a=[\left\{  \pi_{i}\left(  a\right)  \right\} ]_{\omega}.
\]

\end{theorem}

Proof. \ There is an increasing sequence $\left\{  \mathcal{A}_{n}\right\}  $
of finite-dimensional C*-subalgebras of $\mathcal{A}$ whose union
$\mathcal{D}$ is $\left\Vert \cdot\right\Vert _{2}$-dense in $\mathcal{A}$
such that $\mathcal{A}_{1}=\mathbb{C}\cdot1$. Using Lemma
\ref{lemma,lift trace-preserving homo}, for every $i,$ there is a
trace-preserving homomorphism $\pi_{i}:\mathcal{D}\rightarrow\mathcal{M}_{i}$
such that, for every $a\in\mathcal{A},$
\[
a=[\left\{  \pi_{i}\left(  a\right)  \right\} ]_{\omega}.
\]
However, since each $\pi_{i}$ is an isometry in $\left\Vert
\cdot\right\Vert _{2},$ we can extend $\pi_{i}$ uniquely to an
isometry (i.e., trace-preserving) linear map (still called
$\pi_{i}$) from $\mathcal{A}$ to $\mathcal{M}_{i}$. Since
multiplication and the map $x\rightarrow x^{\ast}$ are $\left\Vert
\cdot\right\Vert _{2}$-continuous on bounded sets, it follows that
$\pi_{i}:\mathcal{A}\rightarrow\mathcal{M}_{i}$ is a *-homomorphism,
and that
\[
a=[\left\{  \pi_{i}\left(  a\right)  \right\} ]_{\omega}%
\]
holds for every $a\in\mathcal{A}$.

\end{document}